\documentclass[11pt]{article}

\usepackage{fullpage,amsmath,epsfig}
\usepackage{graphicx}
\usepackage{caption2}

\begin{document}

\newtheorem{definition}{Definition}[section]
\newtheorem{theorem}{Theorem}[section]
\newtheorem{prop}[theorem]{Proposition}
\newtheorem{lemma}[theorem]{Lemma}
\newtheorem{remark}[theorem]{Remark}
\newtheorem{corollary}[theorem]{Corollary}
\newtheorem{observation}[theorem]{Observation}
\newtheorem{claim}[theorem]{Claim}
\newtheorem{conj}[theorem]{Conjecture}

\def\square{\vrule height6pt width7pt depth1pt}
\def\endpf{\hfill\square\bigskip}

\title{Cartesian Products of Regular Graphs are Antimagic}

\author{
\\
  Yongxi~Cheng
\\
\\
\small{Department of Computer Science,
  Tsinghua University, Beijing 100084,
  China}
  \\
[1mm] \small{cyx@mails.tsinghua.edu.cn}
}

\date{} 

\maketitle

\begin{abstract}
An \emph{antimagic labeling} of a finite undirected simple graph with $m$ edges and $n$ vertices is a bijection
from the set of edges to the integers $1,\ldots,m$ such that all $n$ vertex sums are pairwise distinct, where a
vertex sum is the sum of labels of all edges incident with the same vertex. A graph is called \emph{antimagic} if
it has an antimagic labeling. In 1990, Hartsfield and Ringel \cite{HaRi} conjectured that every simple connected
graph, but $K_2$, is antimagic. In this article, we prove that a new class of Cartesian product graphs are
antimagic. In addition, by combining this result and the antimagicness result on toroidal grids (Cartesian
products of two cycles) in \cite{Wan}, all Cartesian products of two or more regular graphs can be proved to be
antimagic.
\\[2mm]
\end{abstract}

\section{Introduction}
All graphs in this paper are finite, undirected and simple. We follow the notation and terminology of \cite{HaRi}.
In 1990, Hartsfield and Ringel \cite{HaRi} introduced the concept of \emph{antimagic} graph. An \emph{antimagic
labeling} of a graph with $m$ edges and $n$ vertices is a bijection from the set of edges to the integers
$1,\ldots,m$ such that all $n$ vertex sums are pairwise distinct, where a vertex sum is the sum of labels of all
edges incident with that vertex. A graph is called antimagic if it has an antimagic labeling. Hartsfield and
Ringel showed that paths $P_n (n\geq 3)$, cycles, wheels, and complete graphs $K_n (n\geq 3)$ are antimagic. They
conjectured that all trees except $K_2$ are antimagic. Moreover, all connected graphs except $K_2$ are antimagic.
These two conjectures are unsettled. In \cite{AKLRY}, Alon et al showed that the latter conjecture is true for all
graphs with $n$ vertices and minimum degree $\Omega (\log n)$. They also proved that complete partite graphs
(other than $K_2$) and $n$-vertex graphs with maximum degree at least $n-2$ are antimagic. In \cite{He}, Hefetz
proved several special cases and variants of the latter conjecture, the main tool used is the Combinatorial
NullStellenSatz (see \cite{Al}). In \cite{Wan}, Wang showed that the toroidal grids, i.e., Cartesian products of
two or more cycles, are antimagic.

In this paper, we prove that the Cartesian products $G_1\times G_2$ of a regular graph $G_1$ and a graph $G_2$ of
bounded degrees are antimagic, provided that the degrees of $G_1$ and $G_2$ satisfy some inequality. By combining
this result and the antimagicness result on the Cartesian products of two cycles in \cite{Wan}, all Cartesian
products of two or more regular graphs (not necessarily connected) can be proved to be antimagic. First, we
introduce another concept about graph labeling called \emph{$\delta$-approximately magic}.

\begin{definition}
\emph{A \emph{$\delta$-approximately magic labeling} of a graph with
$m$ edges is a bijection from the set of edges to the integers
$1,\ldots,m$ such that the difference between the largest and the
smallest vertex sums is at most $\delta$, where a vertex sum is the
sum of labels of all edges incident with that vertex. A graph is
called \emph{$\delta$-approximately magic} if it has a
$\delta$-approximately magic labeling. }
\end{definition}

Thus 0-approximately magic is the same as \emph{magic} in \cite{HaRi}, or \emph{supermagic} in some literature. We
first prove some approximately magicness results on connected regular graphs, the following is proved in Section
2.

\begin{theorem}\label{appromagic}
\textit{If G is an n-vertex k-regular connected graph ($k\geq 1$),
then G is $(\frac {nk} {2}-1)$-approximately magic in case k is odd,
k-approximately magic in case k is even.}
\end{theorem}

Recall that the \emph{Cartesian product} $G_1\times G_2$ of two
graphs $G_1=(V_1, E_1)$ and $G_2=(V_2, E_2)$ is a graph with vertex
set $V_1\times V_2$, and $(u_1,u_2)$ is adjacent to $(v_1,v_2)$ in
$G_1\times G_2$ if and only if $u_1=v_1$ and $u_2v_2\in E_2$, or,
$u_2=v_2$ and $u_1v_1\in E_1$.

Using the approximately magicness results in Theorem
\ref{appromagic}, we prove the following theorem in Section 3.
\begin{theorem}\label{general}
\textit{If $G_1$ is an $n_1$-vertex $k_1$-regular connected graph, and $G_2$ is a graph (not necessarily
connected) with maximum degree at most $k_2$, minimum degree at least one, then $G_1\times G_2$ is antimagic,
provided that $k_1$ is odd and $\frac {k_1^2-k_1}{2}\geq k_2$, or, $k_1$ is even and $\frac {k_1^2}{2}\geq k_2$
and $k_1, k_2$ are not both equal to 2.}
\end{theorem}

By combining Theorem \ref{general} and the antimagicness result on the Cartesian products of two cycles in
\cite{Wan}, the following theorem is obtained in Section 4.

\begin{theorem}\label{main}
\textit{All Cartesian products of two or more regular graphs (not
necessarily connected) are antimagic.}
\end{theorem}

Finally, we give a generalization of Theorems \ref{appromagic} in which $G$ is not necessarily connected, and a
generalization of Theorem \ref{general} in which $G_1$ is not necessarily connected. The following two theorems
are proved in Section 5.

\begin{theorem}\label{uncon-appromagic}
\emph{(generalization of Theorem \ref{appromagic})} \textit{If G is an n-vertex k-regular graph ($k\geq 1$, G is
not necessarily connected), then G is $(\frac {nk} {2}-1)$-approximately magic in case k is odd, $(\frac
{2n}{3}+k-1)$-approximately magic in case k is even.}
\end{theorem}

\begin{theorem}\label{uncon-general}
\emph{(generalization of Theorem \ref{general})} \textit{If $G_1$ is an $n_1$-vertex $k_1$-regular graph, and
$G_2$ is a graph with maximum degree at most $k_2$, minimum degree at least one ($G_1$,$G_2$ are not necessarily
connected), then $G_1\times G_2$ is antimagic, provided that $k_1$ is odd and $\frac {k_1^2-k_1}{2}\geq k_2$, or,
$k_1$ is even and $\frac {k_1^2}{2}>k_2$.}
\end{theorem}

For more results, open problems and conjectures on magic graphs, antimagic graphs and various graph labeling
problems, please see \cite{Ga}.
\\[2mm]
\indent Throughout the paper, we denote by $\lceil x \rceil$
(ceiling of $x$) the least integer that is not less than $x$, denote
by $\lfloor x \rfloor$ (floor of $x$) the largest integer that is
not greater than $x$.
\section{Proof of Theorem \ref{appromagic}}

We begin with some terms and definitions (see \cite{HaRi}). A \emph{walk} in a graph $G$ is an alternating
sequence $v_1e_1v_2e_2\cdots e_{t-1}v_t$ of vertices and edges of $G$, with the property that every edge $e_i$ is
incident with $v_i$ and $v_{i+1}$, for $i=1,\ldots,t-1$. Vertices and edges may be repeated in a walk. A
\emph{trail} in a graph $G$ is a walk in $G$ with the property that no edge is repeated. A \emph{circuit} is a
closed trail, that is a trail whose endpoints are the same vertex. A \emph{cycle} is a circuit with the property
that no vertex is repeated. An \emph{Eulerian circuit} in a graph $G$ is a circuit that contains every edge of
$G$. In order to prove Theorem \ref{appromagic} for the case that $k$ is odd, we need the following theorem
(\cite{HaRi}, pp. 56),

\begin{theorem}\label{listing}
\emph{(part of Listing Theorem)}. \textit{If G is a connected graph
with precisely 2h vertices of odd degree, $h\neq0$, then there exist
h trails in G such that each edge of G is in exactly one of these
trails.}
\end{theorem}

\begin{figure}[t]
\renewcommand{\captionlabelfont}{\bf}
\renewcommand{\captionlabeldelim}{.~}
\centering
\includegraphics[width=145mm]{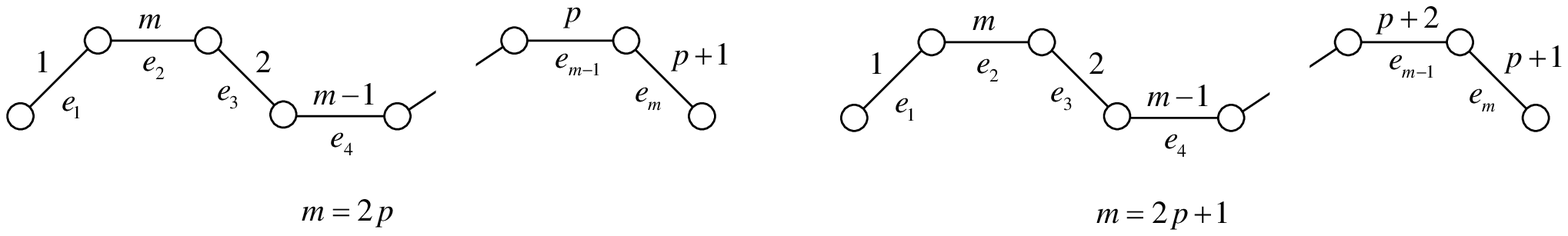}
\renewcommand{\figurename}{Fig.}
\caption{Labeling of the sequence of trails $T : t_1 t_2 \ldots
t_{\frac{n}{2}}$.} \label{fig:trails}
\end{figure}

If $G$ is a connected $n$-vertex regular graph of odd degree $k$, by
Theorem \ref{listing}, there are $n/2$ trails $t_1, t_2, \ldots ,
t_{\frac{n}{2}}$ in $G$, such that each edge of $G$ is in exactly
one of these trails. Denote $|t|$ to be the length (number of edges)
of a trail $t$. Without loss of generality, assume that $|t_1|\geq
|t_2|\geq \ldots \geq |t_{\frac{n}{2}}|$. By concatenating these
trails we get a sequence $T : t_1 t_2 \ldots t_{\frac{n}{2}}$, which
contains all the $m$ $(={\frac{nk}{2}})$ edges of $G$. Number the
edges of $G$ according to their ordering in $T$, let $e_1, e_2,
\ldots, e_m$ be the numbering. Assign the labels $1,2,\ldots, \lceil
\frac{m}{2} \rceil$ to the edges of odd indices $e_1,e_3,\ldots$
etc., and assign the labels $m, m-1,\ldots, \lceil \frac{m}{2}
\rceil+1$ to the edges of even indices $e_2,e_4,\ldots$ etc. (see
Figure \ref{fig:trails}). It is easy to see that for the above
labeling, the sum of any two consecutive edges in $T$ is either
$m+1$ or $m+2$. In addition, if $e$ is the first or the last edge of
a trail, then the largest possible label received by $e$ is at most
$m-\frac{k-1}{2}$ (notice that $|t_1|\geq k$). For each vertex $v$
of $G$, the $k$ edges incident with $v$ can be partitioned into
$\frac{k-1}{2}$ pairs and a singleton, such that each pair is
composed of two consecutive edges within one of the above $n/2$
trails, and the single edge is the first or the last edge of a
trail. Therefore, for the above labeling, the sum received by any
vertex of $G$ is at most $(m-\frac{k-1}{2})+\frac{k-1}{2}\times
(m+2)=m+\frac{k-1}{2}\times (m+1)$, at least $1+\frac{k-1}{2}\times
(m+1)$, implying that this is an $(\frac {nk} {2}-1)$-approximately
magic labeling of $G$. For the case that the degree $k$ is even, we
need the following lemma.

\begin{lemma} \label{lemm:approcycle}
\noindent Every m-vertex connected regular graph of degree 2 (i.e.,
cycle $C_m$) is 2-approximately magic, for $m\geq 3$.
\end{lemma}

\begin{figure}[t]
\renewcommand{\captionlabelfont}{\bf}
\renewcommand{\captionlabeldelim}{.~}
\centering
\includegraphics[width=60mm]{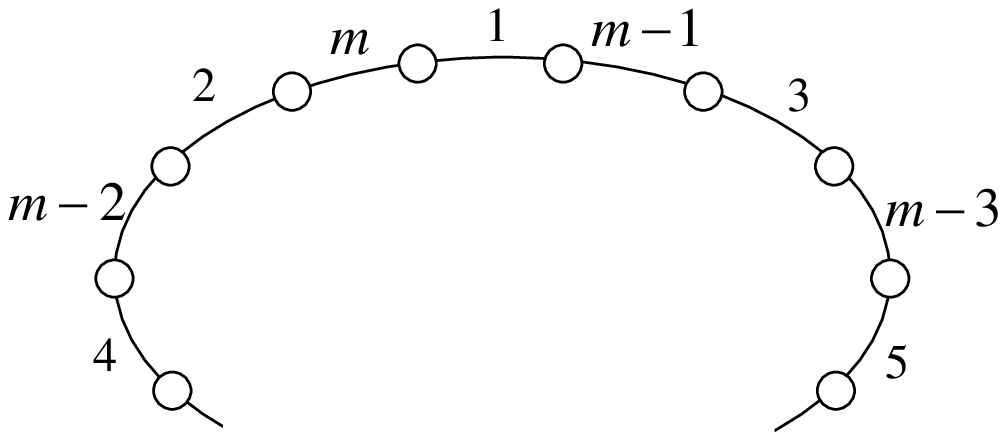}
\renewcommand{\figurename}{Fig.}
\caption{2-Approximately magic labeling of $C_m$} \label{fig:cycles}
\end{figure}

\noindent {\bf Proof:}\, We have the following four cases:
\\[2mm]
\noindent \emph{Case 1.} $m \equiv 1$ (mod 4). Let $m=4t+1$, $t\geq
1$. Partition the labels $1,2,\ldots, m$ into $2t+1$ groups
$(1),(2,3),\ldots,(2t,2t+1),(2t+2,2t+3),\ldots,(m-1,m)$. First
assign label 1 to an arbitrary edge of $C_m$, then assign the labels
$(m, m-1), (2,3), (m-2,m-3), (4,5), \ldots, (2t, 2t+1)$ in a way
that each pair of labels are assigned to the two edges that have
common endpoints with the labeled arc.
\\[2mm]
\noindent \emph{Case 2.} $m \equiv 3$ (mod 4). Let $m=4t+3$, $t\geq
0$. Partition the labels $1,2,\ldots, m$ into $2t+2$ groups
$(1),(2,3),\ldots,(2t,2t+1),(2t+2,2t+3),\ldots,(m-1,m)$. First
assign label 1 to an arbitrary edge of $C_m$, then assign the labels
$(m, m-1), (2,3), (m-2,m-3), (4,5), \ldots, (2t+3, 2t+2)$ in the
same way as in Case 1.
\\[2mm]
\noindent \emph{Case 3.} $m \equiv 0$ (mod 4). Let $m=4t+4$, $t\geq
0$. Partition the labels $1,2,\ldots, m$ into $2t+3$ groups
$(1),(2,3),\ldots,(2t,2t+1),(2t+2),(2t+3,2t+4),\ldots,(m-1,m)$.
First assign label 1 to an arbitrary edge of $C_m$, then assign the
labels $(m, m-1), (2,3), (m-2,m-3), (4,5), \ldots, (2t+4, 2t+3)$ in
the way that each pair of labels are assigned to the two edges that
have common endpoints with the labeled arc, finally assign the label
$2t+2$ to the one non-labeled edge.
\\[2mm]
\noindent \emph{Case 4.} $m \equiv 2$ (mod 4). Let $m=4t+2$, $t\geq
1$. Partition the labels $1,2,\ldots, m$ into $2t+2$ groups
$(1),(2,3),\ldots,(2t,2t+1),(2t+2),(2t+3,2t+4),\ldots,(m-1,m)$.
First assign label 1 to an arbitrary edge of $C_m$, then assign the
labels $(m, m-1), (2,3), (m-2,m-3), (4,5), \ldots, (2t, 2t+1),
(2t+2)$ in the same way as in Case 3.
\\[2mm]
It is easy to see that in any of the above cases, the vertex sums of
$C_m$ are all among $m, m+1,$ and $m+2$, implying the assertion of
the lemma (see Figure \ref{fig:cycles}).
\endpf

Recall that a connected graph with all vertices of even degrees has
an Eulerian circuit. It follows that if $G$ is a connected
$n$-vertex regular graph of even degree $k$, $G$ has an Eulerian
circuit, without loss of generality, say $e_1e_2\ldots e_m$, where
$m={\frac{nk}{2}}$. We label $1,2,\ldots, m$ to this circuit using
the above 2-approximately magic labeling in Lemma
\ref{lemm:approcycle} (here we view this circuit as a cycle). For
each vertex $v$ of $G$, the $k$ edges incident with $v$ can be
partitioned into $k/2$ pairs such that each pair is composed of two
consecutive edges in the Eulerian circuit $e_1e_2\ldots e_m$, thus
the sum of each pair is among $m, m+1,$ and $m+2$. Therefore, for
the above labeling, the sum received by any vertex of $G$ is at
least $\frac{k}{2}\times m$, at most $\frac{k}{2}\times (m+2)$,
implying that this is a $k$-approximately magic labeling of $G$.

\section{Proof of Theorem \ref{general}}

Suppose that $G_1$ is an $n_1$-vertex $k_1$-regular connected graph,
$V(G_1)=\{u_1,u_2,\ldots,u_{n_1}\}$, and $G_2$ is a graph with
maximum degree at most $k_2$, minimum degree at least one,
$V(G_2)=\{v_1,v_2,\ldots,v_{n_2}\}$. Denote by $m_1$ ($=\frac
{k_1n_1}{2}$) and $m_2$ the number of edges of $G_1$ and $G_2$,
respectively.

Let $f: E(G_1\times G_2)\rightarrow \{1,2,\ldots, m_2n_1+m_1n_2\}$
be an edge labeling of $G_1\times G_2$, and denote the induced sum
at vertex $(u,v)$ by $w(u,v)=\sum f((u,v),(y,z))$ , where the sum
runs over all vertices $(y,z)$ adjacent to $(u,v)$ in $G_1\times
G_2$. In the product graph $G_1\times G_2$, at each vertex $(u,v)$,
the edges incident to this vertex can be partitioned into two parts,
one part is contained in a copy of $G_1$ component, and the other
part is contained in a copy of $G_2$ component. Denote by $w_1(u,v)$
and $w_2(u,v)$ the sum at vertex $(u,v)$ restricted to $G_1$
component and $G_2$ component respectively, i.e., $w_1(u,v)=\sum
f((u,v),(y,v))$, where the sum runs over all vertices $y$ adjacent
to $u$ in $G_1$, and $w_2(u,v)=\sum f((u,v),(u,z))$, where the sum
runs over all vertices $z$ adjacent to $v$ in $G_2$. Therefore,
$w(u,v)=w_1(u,v)+w_2(u,v)$.

Given two isomorphic graphs $G$ and $G'$, and two labelings $f$ and $f'$ of $G$ and $G'$ respectively, we call
$f'$ is a \emph{$\delta$-shift} of $f$, if for each edge $e\in E(G)$ and its counterpart $e'\in E(G')$ under the
isomorphism, we have $f'(e')=f(e)+\delta$. Now we will present our labeling of $G_1\times G_2$, which contains two
steps.
\\[2mm]
\noindent \emph{Step 1 (renaming vertices):} Assign labels
$1,2,\ldots,m_1$ to the edges of $G_1$, such that the labeling is
$(\frac {n_1k_1} {2}-1)$-approximately magic if $k_1$ is odd,
$k_1$-approximately magic if $k_1$ is even. Without loss of
generality, we can rename the vertices of $G_1$ such that
$w(u_1)\leq w(u_2)\leq \ldots \leq w(u_{n_1})$, denote this labeling
by $L_1$. Assign labels $1,n_1+1,2n_1+1,\ldots,(m_2-1)n_1+1$
arbitrarily to the edges of $G_2$. Similarly, rename the vertices of
$G_2$ such that $w(v_1)\leq w(v_2)\leq \ldots \leq w(v_{n_2})$,
denote this labeling by $L_2$.
\\[2mm]
\noindent \emph{Step 2 (labeling on $G_1\times G_2$):} Assign labels $m_2n_1+1,m_2n_1+2,\ldots,m_2n_1+m_1n_2$ to
the edges that are contained in copies of $G_1$ component. For the $i$-th $G_1$ component (with vertices
$(u_1,v_i)$, $(u_2,v_i)$,\ldots, $(u_{n_1},v_i)$), label its edges with
$m_2n_1+(i-1)m_1+1,m_2n_1+(i-1)m_1+2,\ldots,m_2n_1+(i-1)m_1+m_1$, such that the labeling is an
$[m_2n_1+(i-1)m_1]$-shift of $L_1$, under the natural isomorphism, for $i=1,\ldots,n_2$. Since $G_1$ is regular,
we have $w_1(u_1,v_i)\leq w_1(u_2,v_i)\leq \ldots \leq w_1(u_{n_1},v_i)$, for $i=1,\ldots,n_2$.

Assign labels $1,2,\ldots,m_2n_1$ to the edges that are contained in copies of $G_2$ component. For the $j$-th
$G_2$ component (with vertices $(u_j,v_1)$, $(u_j,v_2)$,\ldots, $(u_j,v_{n_2})$), label its edges with
$j,n_1+j,2n_1+j,\ldots,(m_2-1)n_1+j$, such that the labeling is a $(j-1)$-shift of $L_2$, under the natural
isomorphism, for $j=1,\ldots,n_1$. From the way we name the vertices of $G_2$, we have $w_2(u_1,v_1)\leq
w_2(u_1,v_2)\leq \ldots \leq w_2(u_1,v_{n_2})$.\\

In what follows we will prove that for the above labeling, if $k_1$
is odd and $\frac {k_1^2-k_1}{2}\geq k_2$, or, if $k_1$ is even and
$\frac {k_1^2}{2}\geq k_2$ and $k_1, k_2$ are not both equal to 2,
then
\begin{eqnarray} \label{monotone}
\nonumber
&&w(u_1,v_1)<w(u_2,v_1)<\ldots\ldots\ldots\ldots<w(u_{n_1},v_1)<
\\
&&w(u_1,v_2)<w(u_2,v_2)<\ldots\ldots\ldots\ldots<w(u_{n_1},v_2)<
\\
&&\;\;\;\;\;\; \nonumber
\ldots\ldots\ldots\ldots\ldots\ldots\ldots\ldots\ldots\ldots\ldots\ldots\ldots\ldots\ldots
\\
\nonumber
&&w(u_1,v_{n_2})<w(u_2,v_{n_2})<\ldots\ldots\ldots\ldots<w(u_{n_1},v_{n_2}),
\end{eqnarray}
implying that the above labeling is antimagic.\\

For each $i\in \{1,\ldots,n_2\}$, we have $w_1(u_1,v_i)\leq
w_1(u_2,v_i)\leq \ldots \leq w_1(u_{n_1},v_i)$, and
$w_2(u_1,v_i)<w_2(u_2,v_i)<\ldots <w_2(u_{n_1},v_i)$ since
$w_2(u_{j+1},v_i)-w_2(u_j,v_i)=d(v_i)$, where $d(v_i)\geq 1$ is the
degree of $v_i$ in $G_2$, $j=1,\ldots,n_1-1$. It follows that
$w(u_1,v_i)<w(u_2,v_i)<\ldots<w(u_{n_1},v_i)$, for $i=1,\ldots,n_2$.
In order to prove $w(u_1,v_{i+1})>w(u_{n_1},v_i)$, for
$i=1,\ldots,n_2-1$, there are two cases.
\\[2mm]
\emph{Case 1.} $k_1$ is odd. For each $i\in \{1,\ldots,n_2-1\}$, we
have $w(u_1,v_{i+1})\geq w(u_1,v_i)+\frac{n_1k_1^2}{2}$ since
$w_1(u_1,v_{i+1})=w_1(u_1,v_i)+m_1k_1=w_1(u_1,v_i)+\frac{n_1k_1^2}{2}$
(notice that the labeling on the $(i+1)$-th $G_1$ component is an
$m_1$-shift of the labeling on the $i$-th $G_1$ component) and
$w_2(u_1,v_{i+1})\geq w_2(u_1,v_i)$. In addition, we have
$w(u_{n_1},v_i)\leq w(u_1,v_i)+(\frac{n_1k_1}{2}-1)+k_2(n_1-1)$
since $w_1(u_{n_1},v_i)\leq w_1(u_1,v_i)+(\frac{n_1k_1}{2}-1)$
(notice that $G_1$ is regular and $L_1$ is $(\frac {n_1k_1}
{2}-1)$-approximately magic when $k_1$ is odd), and
$w_2(u_{n_1},v_i)=w_2(u_1,v_i)+d(v_i)(n_1-1)\leq
w_2(u_1,v_i)+k_2(n_1-1)$. It follows that
$w(u_1,v_{i+1})-w(u_{n_1},v_i)\geq
(w(u_1,v_i)+\frac{n_1k_1^2}{2})-(w(u_1,v_i)+(\frac{n_1k_1}{2}-1)+k_2(n_1-1))=n_1(\frac
{k_1^2-k_1}{2}-k_2)+1+k_2>0$, for $i=1,\ldots,n_2-1$.
\\[2mm]
\emph{Case 2.} $k_1$ is even. Similarly, for each $i\in
\{1,\ldots,n_2-1\}$, we have $w(u_1,v_{i+1})\geq
w(u_1,v_i)+\frac{n_1k_1^2}{2}$ since
$w_1(u_1,v_{i+1})=w_1(u_1,v_i)+m_1k_1=w_1(u_1,v_i)+\frac{n_1k_1^2}{2}$
and $w_2(u_1,v_{i+1})\geq w_2(u_1,v_i)$. In addition,
$w(u_{n_1},v_i)\leq w(u_1,v_i)+k_1+k_2(n_1-1)$ holds since
$w_1(u_{n_1},v_i)\leq w_1(u_1,v_i)+k_1$ ($L_1$ is
$k_1$-approximately magic when $k_1$ is even) and
$w_2(u_{n_1},v_i)=w_2(u_1,v_i)+d(v_i)(n_1-1)\leq
w_2(u_1,v_i)+k_2(n_1-1)$. It follows that
$w(u_1,v_{i+1})-w(u_{n_1},v_i)\geq
(w(u_1,v_i)+\frac{n_1k_1^2}{2})-(w(u_1,v_i)+k_1+k_2(n_1-1))=n_1(\frac
{k_1^2}{2}-k_2)+k_2-k_1$.

If $\frac {k_1^2}{2}>k_2$, since $k_1$ is even, $\frac
{k_1^2}{2}-k_2\geq 1$, then $w(u_1,v_{i+1})-w(u_{n_1},v_i)\geq
n_1(\frac {k_1^2}{2}-k_2)+k_2-k_1\geq n_1+k_2-k_1>0$ (since
$n_1>k_1$). If $\frac {k_1^2}{2}=k_2$, since $k_1,k_2$ are not both
equal to 2, we have $k_1>2$, it follows that
$w(u_1,v_{i+1})-w(u_{n_1},v_i)\geq k_2-k_1=\frac {k_1^2}{2}-k_1>0$.
Thus, in any case, we have $w(u_1,v_{i+1})-w(u_{n_1},v_i)>0$, for
$i=1,\ldots,n_2-1$.
\\[2mm]

Therefore, (\ref{monotone}) holds, implying the assertion of Theorem
\ref{general}.

\section{Proof of Theorem \ref{main}}

Since the Cartesian product preserves regularity, we only need to prove that all Cartesian products of two regular
graphs are antimagic. We first prove Theorem \ref{main} for the case that $G_1$ and $G_2$ are both connected, then
we generalize the proof to the case where $G_1$ and $G_2$ are not necessarily connected.

\subsection{Connected Case}
Suppose that $G_1$ is an $n_1$-vertex $k_1$-regular connected graph, and $G_2$ is an $n_2$-vertex $k_2$-regular
connected graph. Without loss of generality, assume that $k_1\geq k_2$. Furthermore, we may assume $k_1\geq 2$
since $K_2\times K_2$ can be easily verified as antimagic.  If $k_1=2$ and $k_2=1$, by Theorem \ref{general},
$G_1\times G_2$ is antimagic. If $k_1=2$ and $k_2=2$, then $G_1\times G_2$ is a toroidal grid graph and its
antimagicness is proved in \cite{Wan}. For $k_1\geq 3$, if $k_1$ is odd, then $\frac {k_1^2-k_1}{2}\geq k_1\geq
k_2$; if $k_1$ is even, then $k_1\geq 4$, $\frac {k_1^2}{2}>k_1\geq k_2$. Thus by Theorem \ref{general},
$G_1\times G_2$ is antimagic.

\subsection{Unconnected Case}
Denote by $c_1$ and $c_2$ the numbers of connected components of
$G_1$ and $G_2$, respectively. It is easy to see that the number of
connected components of $G_1\times G_2$ is $c=c_1\times c_2$, and
each of its connected components is a $(k_1+k_2)$-regular graph
(which is product of one $k_1$-regular connected graph and one
$k_2$-regular connected graph). Let $m_1,m_2,\ldots,m_c$ be the
numbers of edges of these connected components $C_1,C_2,\ldots,C_c$.
The labeling of $G_1\times G_2$ goes as follows. Assign
$1,2,\ldots,m_1$ to the edges of $C_1$, assign
$m_1+1,m_1+2,\ldots,m_1+m_2$ to the edges of $C_2$, \ldots\ldots,
and assign $m_1+\ldots +m_{c-1}+1, m_1+\ldots
+m_{c-1}+2,\ldots,m_1+\ldots +m_{c-1}+m_c$ to the edges of $C_c$,
such that the labeling of each connected component is antimagic
(this can be achieved because of the previous proof for the case
where $G_1$ and $G_2$ are both connected and the regularity of each
component). The whole labeling of $G_1\times G_2$ is antimagic,
since between any two different components, any sum of $k_1+k_2$
labels from a group of larger labels must be greater than any sum of
$k_1+k_2$ labels from a group of smaller labels.

\section{Generalizations of Theorem \ref{appromagic} and \ref{general}}

In this section, we will prove Theorems \ref{uncon-appromagic}, a generalization of Theorem \ref{appromagic} in
which $G$ is not necessarily connected, and Theorem \ref{uncon-general}, a generalization of Theorem \ref{general}
in which $G_1$ is not necessarily connected.

\subsection{Proof of Theorem \ref{uncon-appromagic}}
For the case $k$ is odd, by Theorem \ref{listing} (Listing), for
each connected component of $G$ (which is a connected $k$-regular
graph), if it has $n_i$ vertices, we can decompose it into
$\frac{n_i}{2}$ trails. By running this decomposition over all
connected components of $G$, we can get a total number of
$\frac{n}{2}$ trails, such that each edge of $G$ is in exactly one
of these trails. It is easy to see that the largest length of these
trails is at least $k$. We concatenate these trails into a sequence
in the ordering of nonincreasing lengths, and label the sequence in
the same way as in Theorem \ref{appromagic}, which results in an
$(\frac {nk} {2}-1)$-approximately magic labeling of $G$. For the
case $k$ is even, we first prove the following lemma.

\begin{lemma} \label{lemm:multicycles}
\noindent If G is an n-vertex graph consisting of vertex-disjoint
cycles of odd sizes (numbers of edges), then G is $\lceil
\frac{2n}{3} \rceil$-approximately magic, for $n\geq 3$.
\end{lemma}

\noindent {\bf Proof:}\, Suppose that $G$ is composed of $l$ cycles
$C_1$,$C_2$,\ldots,$C_l$ (of sizes $n_1,n_2,\ldots, n_l$, where
$n_1\geq n_2\geq\ldots \geq n_l\geq 3$ are odd numbers, and
$n_1+\cdots+n_l=n$). Let $n=3t+\varepsilon$, $t\geq1$, $\varepsilon
\in \{0,1,2\}$. We partition the labels $1,\ldots,n$ into three
groups $1,2,\ldots,t$ and $t+1,\ldots,2t+\varepsilon$ and
$2t+\varepsilon+1,2t+\varepsilon+2,\ldots,3t+\varepsilon$. Let $A:
a_1,a_2,\ldots,a_t$ denote the sequence $1,2,\ldots,t$; let $B:
b_1,b_2,\ldots,b_{t+\varepsilon}$ denote the sequence
$2t+\varepsilon,2t+\varepsilon-1,\ldots,t+1$; and let $C:
c_1,c_2,\ldots,c_t$ denote the sequence
$2t+\varepsilon+1,2t+\varepsilon+2,\ldots,3t+\varepsilon$. It is
easy to see that $2t+\varepsilon+2 \leq a_i+c_j\leq 4t+\varepsilon$,
$a_i+b_i=2t+\varepsilon+1$, and $b_i+c_i=4t+2\varepsilon+1$, for
$i,j=1,2,\ldots,t$. In addition, $2t+3 \leq b_i+b_j\leq
4t+2\varepsilon-1$, for $i\neq j$, $i,j=1,2,\ldots,t$.
\\[2mm]
\indent Let $m_i=\frac{n_i-1}{2}$, $i=1,2,\ldots,t$. We will present
a labeling on $G$, which goes as follows. Label the cycles
$C_1$,$C_2$,\ldots,$C_l$ one by one. For the $i$-th cycle $C_i$,
pick the $m_i$ smallest elements from the current (remained)
$A$-sequence and the $m_i$ smallest elements from the current
(remained) $C$-sequence, if at this moment there are at least $m_i$
elements remained in $A$ (also $C$). Otherwise, pick all the
remained elements of the two sequences. Specifically, we have the
following two cases.
\\[2mm]
\noindent \emph{Case 1.} At the beginning of the labeling of $C_i$,
there are at least $m_i$ elements remained in the current $A$ (also
$C$) sequence. Denote by $a_{s_i+1},a_{s_i+2},\ldots,a_{s_i+m_i}$
and $c_{s_i+1},c_{s_i+2},\ldots,c_{s_i+m_i}$ (where $s_1=0$, and
$s_i=m_1+\cdots+m_{i-1}$ for $1<i\leq l$) the $m_i$ smallest
elements of the current $A$ (and $C$) sequence. Pick $b_{s_i+m_i}$
from the current $B$-sequence, and label the edges of $C_i$
sequentially with $b_{s_i+m_i}$, $c_{s_i+1}$, $a_{s_i+1}$,
$c_{s_i+2}$, $a_{s_i+2}$,\ldots,$c_{s_i+m_i}$, $a_{s_i+m_i}$, then
remove these elements from their sequences. Since
$3t+\varepsilon+2\leq b_{s_i+m_i}+c_{s_i+1}\leq 4t+2\varepsilon+1$,
for the above labeling, each vertex sum of $C_i$ is at least
$2t+\varepsilon+1$, and at most $4t+2\varepsilon+1$.
\\[2mm]
\noindent \emph{Case 2.} At the beginning of the labeling of $C_i$,
the number of elements remained in the current $A$ (also $C$)
sequence is less than $m_i$. In this case we must have $n_1\geq 5$
(otherwise all cycles are `triangles', i.e. consisting of 3 edges,
in our labeling each triangle will be labeled by three elements, and
exactly one element from each sequence, which is a contradiction).
Without loss of generality, we can assume that $l\geq 2$, since if
$l=1$, $G$ has been proved to be $2$-approximately magic in Lemma
\ref{lemm:approcycle}.

If the current $A$ (also $C$) sequence is empty, then label the
remained non-labeled cycles arbitrarily using elements remained in
$B$-sequence. Otherwise, pick all the elements
$a_{s_i+1},a_{s_i+2},\ldots,a_t$ and
$c_{s_i+1},c_{s_i+2},\ldots,c_t$ from the current $A$ and $C$
sequences. At this moment, besides $b_t$ (where $t\geq 2$ since
$l\geq 2$), $b_1$ is unused (if $i=1$, since $t\geq 2$, we have
$b_1$ distinct from $b_t$ and unused; if $i>1$, since $n_1\geq 5$,
$b_1$ has not been used for labeling $C_1$, thus is unused). Remove
$b_t$ and $b_1$ from the current $B$-sequence, and label the
elements $b_t$, $c_{s_i+1}$, $a_{s_i+1}$, $c_{s_i+2}$,
$a_{s_i+2}$,\ldots,$c_t$, $a_t$, $b_1$ sequentially to an arc of
consecutive edges of $C_i$. Then, label the remained non-labeled
edges of $C_i$ using arbitrary elements remained in $B$-sequence,
and remove these elements from $B$. Since $3t+\varepsilon+2\leq
b_t+c_{s_i+1}\leq 4t+2\varepsilon+1$, and $a_t+b_1=3t+\varepsilon$,
we have that for the above labeling, each vertex sum of $C_i$ is at
least $2t+\varepsilon+1$, and at most $4t+2\varepsilon+1$.
\\[2mm]
\indent Therefore, for the above labeling, the vertex sums of $G$
are at least $2t+\varepsilon+1$ (which is $\lceil \frac{2n}{3}
\rceil+1$), at most $4t+2\varepsilon+1$ (which is $2\lceil
\frac{2n}{3} \rceil+1$), implying that the differences between
vertex sums of $G$ are at most $\lceil \frac{2n}{3} \rceil$. \endpf

\begin{remark}
\noindent $\lceil \frac{2n}{3} \rceil$ obtained in Lemma \ref{lemm:multicycles} is actually asymptotically best
possible. Consider the case that $G$ is consisting of $\frac{n}{3}$ `triangles'. Suppose that label 1 is assigned
to an edge $v_1v_2$ of a triangle $v_1v_2v_3$, if the edge $v_2v_3$ or $v_1v_3$ is assigned with a label
$l>\frac{2n}{3}$, then the difference of the two vertex sums of $v_3$ and $v_1$, or $v_3$ and $v_2$ will be at
least $\frac{2n}{3}$. Similarly, suppose that label $n$ is assigned to an edge $v_4v_5$ of a triangle $v_4v_5v_6$,
if the edge $v_4v_6$ or $v_5v_6$ is assigned with a label $l\leq \frac{n}{3}$, then the difference of the two
vertex sums of $v_5$ and $v_6$, or $v_4$ and $v_6$ will be at least $\frac{2n}{3}$. If neither of the above two
cases happens, then the vertex sum of $v_1$ or $v_2$ is at most $\frac{2n}{3}$, and the vertex sum of $v_4$ is at
least $\frac{4n}{3}$, thus, the difference of the two vertex sums of $v_4$ and $v_1$, or $v_4$ and $v_2$ is at
least $\frac{2n}{3}$.
\end{remark}

Now we will prove Theorem \ref{uncon-appromagic} for the case that
$k$ is even. Since $k$ is even, $G$ is an \emph{even graph} (a graph
with all vertices having even degrees), it follows that $G$ can be
decomposed into edge-disjoint simple cycles. In addition, two cycles
having a common vertex can be merged into one circuit. Therefore, by
repeating merging two cycles of odd sizes that having a common
vertex into an even circuit, finally we will obtain a collection of
$s$ $(\geq 0)$ even circuits $P_1$, $P_2$, \ldots, $P_s$ (of sizes
$2m_1,2m_2,\ldots,2m_s$), together with a collection of $t$ ($\geq
0$) vertex-disjoint odd cycles $Q_1$, $Q_2$, \ldots, $Q_t$ (of sizes
$n_1,n_2,\ldots, n_t$, and $n_1+n_2+\cdots +n_t\leq n$), such that
each edge of $G$ is in exactly one of these circuits or cycles.

Let $m=\frac{nk}{2}$ be the number of edges of $G$. First we label
the even circuits $P_1$, $P_2$, \ldots, $P_s$. By viewing these
circuits as cycles, using the 2-approximately magic labeling in
Lemma \ref{lemm:approcycle}, we assign labels $1,2,\ldots,m_1$ and
$m,m-1,\ldots,m-m_1+1$ to $P_1$, assign labels
$m_1+1,m_1+2,\ldots,m_1+m_2$ and $m-m_1,m-m_1-1,\ldots,m-m_1-m_2+1$
to $P_2$, \ldots\ldots, and assign labels $m_1+\ldots+m_{s-1}+1$, \
$m_1+\ldots+m_{s-1}+2$, \ldots\ldots, $m_1+\ldots+m_{s-1}+m_s$ and
$m-m_1-\ldots-m_{s-1}$,
 \ $m-m_1-\ldots-m_{s-1}-1$, \ldots\ldots, $m-m_1-\ldots-m_{s-1}-m_s+1$ to
$P_s$. Thus, the sum of any two consecutive edges of circuit $P_i$
($i=1,\ldots,s$) is among $m$, $m+1$, and $m+2$.

Let $m^{*}=m_1+m_2+\ldots+m_s$, and $n^{*}=n_1+n_2+\ldots+n_t$. If $n^{*}=0$ (i.e., there is no odd cycle),
similarly as in Theorem \ref{appromagic}, the above labeling of $G$ can be proved to be $k$-approximately magic,
by partitioning the $k$ edges incident with any vertex of $G$ into $k/2$ pairs such that each pair is composed of
two consecutive edges in some circuit $P_i$ ($i\in \{1,\ldots,s\}$). Otherwise, we have $n^{*}\geq 3$. Assign the
remained labels $m^{*}+1, m^{*}+2, \ldots, m^{*}+n^{*}$ to the vertex-disjoint odd cycles $Q_1$, $Q_2$, \ldots,
$Q_t$, using the $\lceil \frac{2n^{*}}{3} \rceil$-approximately magic labeling in Lemma \ref{lemm:multicycles}.
Since $2m^{*}+n^{*}=m$, and $\lfloor\frac{n^{*}}{3}\rfloor+\lceil\frac{2n^{*}}{3}\rceil=n^{*}$ for all integers
$n^{*}\geq 1$, it follows that the sum of any two consecutive edges of these odd cycles is at least
$2m^{*}+\lceil\frac{2n^{*}}{3}\rceil+1=m+1-\lfloor\frac{n^{*}}{3}\rfloor$ $(\leq m)$, and at most $2m^{*}+2\lceil
\frac{2n^{*}}{3} \rceil+1=m+1-\lfloor\frac{n^{*}}{3}\rfloor+\lceil\frac{2n^{*}}{3}\rceil$ $(\geq m+2)$. Therefore,
for the whole labeling of $G$, the sum received by any vertex of $G$ is at least $m\times
\frac{k-2}{2}+(m+1-\lfloor\frac{n^{*}}{3}\rfloor)$, at most $(m+2)\times
\frac{k-2}{2}+(m+1-\lfloor\frac{n^{*}}{3}\rfloor+\lceil \frac{2n^{*}}{3} \rceil)$. Since $n^{*}\leq n$, the whole
labeling of $G$ is $(\frac {2n}{3}+k-1)$-approximately magic.

\subsection{Proof of Theorem \ref{uncon-general}}

If $k_1=2$, since $\frac {k_1^2}{2}>k_2$, $k_2=1$, $G_2$ is
1-regular, by Theorem \ref{main}, $G_1\times G_2$ is antimagic. In
what follows we assume that $k_1\geq 3$.
\\[2mm]
\indent We do the same labeling on $G_1\times G_2$ as in  Theorem
\ref{general} (when $k_1$ is even, the labeling $L_1$ on $G_1$ here
is $(\frac {2n_1}{3}+k_1-1)$-approximately magic). We will prove
that for this labeling, (\ref{monotone}) still holds if $k_1\geq 3$
is odd and $\frac {k_1^2-k_1}{2}\geq k_2$, or, if $k_1\geq 4$ is
even and $\frac {k_1^2}{2}>k_2$.
\\[2mm]
\indent $w(u_1,v_i)<w(u_2,v_i)<\ldots<w(u_{n_1},v_i)$ can be proved by using the same argument in Theorem
\ref{general}, for $i=1,\ldots,n_2$ . In order to prove $w(u_1,v_{i+1})-w(u_{n_1},v_i)>0$, for $i=1,\ldots,n_2-1$,
there are two cases.
\\[2mm]
\emph{Case 1.} $k_1$ is odd. Since $G_1$ is still $(\frac {n_1k_1} {2}-1)$-approximately magic, by using the same
argument in Theorem \ref{general}, we can obtain that $w(u_1,v_{i+1})-w(u_{n_1},v_i)>0$, for $i=1,\ldots,n_2-1$.
\\[2mm]
\emph{Case 2.} $k_1$ is even (thus $k_1 \geq 4$). $G_1$ is $(\frac
{2n_1}{3}+k_1-1)$-approximately magic. For each $i\in
\{1,\ldots,n_2-1\}$, we have $w(u_1,v_{i+1})\geq
w(u_1,v_i)+\frac{n_1k_1^2}{2}$ since
$w_1(u_1,v_{i+1})=w_1(u_1,v_i)+\frac{n_1k_1^2}{2}$ and
$w_2(u_1,v_{i+1})\geq w_2(u_1,v_i)$. In addition,
$w(u_{n_1},v_i)\leq w(u_1,v_i)+(\frac {2n_1}{3}+k_1-1)+k_2(n_1-1)$
since $w_1(u_{n_1},v_i)\leq w_1(u_1,v_i)+(\frac {2n_1}{3}+k_1-1)$
and $w_2(u_{n_1},v_i)=w_2(u_1,v_i)+d(v_i)(n_1-1)\leq
w_2(u_1,v_i)+k_2(n_1-1)$. Therefore,
$w(u_1,v_{i+1})-w(u_{n_1},v_i)\geq
(w(u_1,v_i)+\frac{n_1k_1^2}{2})-(w(u_1,v_i)+(\frac
{2n_1}{3}+k_1-1)+k_2(n_1-1))=n_1(\frac {k_1^2}{2}-\frac
{2}{3}-k_2)+k_2-k_1+1$.

Since $k_2<\frac {k_1^2}{2}$, there are two cases: $k_2\leq \frac
{k_1^2}{2}-2$ or $k_2=\frac {k_1^2}{2}-1$. If $k_2\leq \frac
{k_1^2}{2}-2$, $w(u_1,v_{i+1})-w(u_{n_1},v_i)\geq n_1(\frac
{k_1^2}{2}-\frac {2}{3}-k_2)+k_2-k_1+1>n_1+k_2-k_1>0$ (since
$n_1>k_1$). If $k_2=\frac {k_1^2}{2}-1$,
$w(u_1,v_{i+1})-w(u_{n_1},v_i)\geq n_1(\frac {k_1^2}{2}-\frac
{2}{3}-k_2)+k_2-k_1+1>\frac {k_1^2}{2}-k_1>0$ (since $k_1\geq 4$).
Thus, in either case, we have
$w(u_1,v_{i+1})-w(u_{n_1},v_i)>0$, for $i=1,\ldots,n_2-1$.\\

Therefore, (\ref{monotone}) holds, the labeling for $k_1\geq 3$ is
antimagic.

\section{Concluding Remarks and Open Problems}

Since the Eulerian circuit of an Eulerian graph (consequently the trails in the Listing Theorem) can be
efficiently computed, the proofs in this paper provide efficient algorithms for finding the antimagic labelings.

It is easy to see that, for cycles, the $2$-approximately magicness result in Lemma \ref{lemm:approcycle} is best
possible (i.e., 2 can not be improved to 0 or 1). For $n$-vertex $k$-regular ($k>2$) connected graphs, it may be
interesting to prove that they are $\delta$-approximately magic, where $\delta<(\frac {nk} {2}-1)$ in case $k$ is
odd, or $\delta<k$ in case $k$ is even, or, to prove some lower bounds on $\delta$.

\end{document}